\documentclass[preprint,10pt]{elsarticle}
\usepackage{amsmath,amssymb,amsthm,mathrsfs,dsfont}
\usepackage[margin=3cm]{geometry}
\usepackage{titlesec,hyperref}
\usepackage{color}
\usepackage{fancyhdr}
\titleformat{\subsection}{\it}{\thesubsection.\enspace}{1pt}{}

\makeatletter
\def\ps@pprintTitle{%
   \let\@oddhead\@empty
   \let\@evenhead\@empty
   \let\@oddfoot\@empty
   \let\@evenfoot\@oddfoot
}
\makeatother

\newtheorem{theo}{Theorem}[section]
\newtheorem{lemm}[theo]{Lemma}

\newtheorem{rema}[theo]{Remark}
\numberwithin{equation}{section}
\def\bel{\begin{equation}\label}

\def\eeq{\end{equation}}

\newcommand\N{{\mathbb N}\,}
\newcommand\R{{\mathbb R}\,}
\newcommand{\al}{{\alpha}}

\newcommand{\na}{{\nabla}}
\newcommand{\pa}{{\partial}}
\newcommand\tr{{\rm tr}\,}
\newcommand{\beq}{\begin{equation}}
\newcommand{\beno}{\begin{equation*}}
\newcommand{\eeno}{\end{equation*}}

\let\pa=\partial
\let\al=\alpha

\let\d=\delta

\let\lam=\lambda

\let\f=\frac

\let\Om=\Omega

\let\tri=\triangle
\let\ep=\epsilon
\def\na{\nabla}
\def\dive{\mathop{\rm div}\nolimits}
\def\exp{\mathop{\rm exp}\nolimits}
\def\det{\mathop{\rm det}\nolimits}

\def\tr{\mathop{\rm tr}\nolimits}

\allowdisplaybreaks

\begin{document}
\begin{frontmatter}
\title{Blow up and global existence for the periodic Phan-Thein-Tanner model}
\author[add1]{Yuhui Chen}\ead{chenyh339@mail.sysu.edu.cn}
\author[add2] {Wei Luo\corref{cor1}}\ead{luowei23@mail2.sysu.edu.cn}\cortext[cor1]{Corresponding author.}
\author[add2]{Zheng-an Yao}\ead{mcsyao@mail.sysu.edu.cn}
\address[add1]{School of  Aeronautics and Astronautic, Sun Yat-sen University, Guangzhou, 510275, CHINA}
\address[add2]{School of Mathematics, Sun Yat-sen University, Guangzhou, 510275, CHINA}

\begin{abstract}
In this paper, we mainly investigate the Cauchy problem for the periodic Phan-Thein-Tanner (PTT) model. This model is derived from network theory for the polymeric fluid. We prove that the strong solutions of PTT model will blow up in finite time if the trace of initial stress tensor $\tr \tau_0(x)$ is negative. This is a very different with the other viscoelastic model. On the other hand, we obtain the global existence result with small initial data when $\tr \tau_0(x)\geq c_0>0$ for some $c_0$. Moreover, we study about the large time behaviour.

\noindent {\it 2010 AMS Classification}: 35A01, 35B45, 35Q35, 76A05, 76D03.
\end{abstract}


\end{frontmatter}
\vspace*{10pt}

\section{Introduction}
In this paper, we consider the initial value problem for the following periodic incompressible Phan-Thein-Tanner (PTT) model\cite{Thien1977,Thien1978}:
\begin{align}\label{PTT}\tag{PTT}
\left\{
\begin{array}{ll}
u_{t}+u\cdot \na u-\mu \tri u+\na p=\mu_1\dive \tau, \quad (t, x) \in \mathbb{R}^{+}\times \mathbb{T}^3, \\[1ex]
\tau_t+u\cdot \na\tau+ (a+b\tr\tau)\tau+Q(\tau,\na u)=\mu_2 D(u),  \quad (t, x) \in \mathbb{R}^{+}\times \mathbb{T}^3, \\[1ex]
\dive u=0,  \quad (t, x) \in \mathbb{R}^{+}\times \mathbb{T}^3, \\[1ex]
u|_{t=0}=u_{0}(x),    \tau|_{t=0}=\tau_0(x), \quad x \in  \mathbb{T}^3. \\[1ex]
\end{array}
\right.
\end{align}
Here $u$ stands for the velocity and $p$ is the scalar pressure of fluid, $\tau$ is the stress tensor. $D(u)$ is the symmetric part of $\na u$, that is
\[D(u)=\f12(\na u+(\na u)^T).\]
$Q(\tau,\na u)$ is a given bilinear form
\[Q(\tau,\na u)=\tau \Om(u)-\Om(u)\tau+\lam(D(u)\tau+\tau D(u)),\]
where $\Om(u)$ is the skew-symmetric part of $\na u$, namely
\[\Om(u)=\f12(\na u-(\na u)^T).\]
$\mu>0$ is the viscosity coefficient and $\mu_1$ is the elastic coefficient. $a$ and $\mu_2$ are associated to the Debroah number $De=\frac{\mu_2}{a}$, which indicates the relation between the characteristic flow time and elastic time\cite{Bird1977}. $\lam\in[-1,1]$ is a physical parameter. In particular, we call the system co-rotational case when $\lam=0$. $b\geq 0$ is a constant relate to the rate of creation or destruction for the polymeric network junctions.

If $b=0$. the system \eqref{PTT} reduce to the famous Oldroyd-B model (See \cite{Oldroyd}) which has been studied widely. Let us review some mathematical results for the related Oldroyd type model. C. Guillop\'e and J. C. Saut \cite{Guillope1990-NA,Guillope1990} proved the existence of local strong solutions and the global existence of one dimensional shear flows. In \cite{Fernandez-Cara}, E. Fern\'andez-Cara, F. Guill\'en and R. Ortega studied the local well-posedness in Sobolev spaces. J. Chemin and N. Masmoudi \cite{Chemin2001} proved the local well-posedness in critical Besov spaces and give a low bound for the lifespan. In the co-rotational case, P. L. Lions and N. Masmoudi \cite{Lions-Masmoudi} proved the global existence of weak solutions. In \cite{Lin-Liu-Zhang2005}, F. Lin, C. Liu and P. Zhang proved that if the initial data is a small perturbation around equilibrium, then the strong solution is global in time. The similar results were obtained in several papers by virtue of different methods, see Z. Lei and Y. Zhou \cite{Lei-Zhou2005}, Z. Lei, C. Liu and Y. Zhou \cite{Lei2008}, T. Zhang and D. Fang \cite{Zhang-Fang2012}, Y. Zhu \cite{Zhu2018}. D. Fang, M.Hieber and R. Zi proved the global existence of strong solutions with a class of large data \cite{Fang-Zi2013,Fang-Zi2016}. For the Oldroyd-B model, the global existence of strong solutions in two dimension without small conditions is still an open problem.

In this paper, we focus on the PTT model $b\neq 0$. To our knowledge, there are a lot of numerical results about the PTT model (See, \cite{Oliveira1999,Mu-Zhao2012,Mu-Zhao2013,Bautista,Tamaddon-Jahromi}). However, there is no any well-posedness results about the PTT model. The nonlinear term $(\tr\tau) \tau$ in the PTT model will leads to some interesting phenomenon that is quiet different between the Oldroyd-B model. By virtue of the characteristic method, we prove that the strong solution of \eqref{PTT} will blow up in finite time when the initial data $\tr \tau_0<0$. This phenomenon can not be founded in other viscoelastic model.

On the other hand, when $\tr \tau_0$ has a positive low bound $c_0$, we can prove the global existence of strong solution with small initial data. The idea is inspired by the method applied in \cite{Zhu2018}. The main different is to deal with the nonlinear term $(\tr\tau) \tau$. In \cite{Zhu2018}, the author observe that the linearized system of the Oldroyd-B model has the form
\beno
\left\{
\begin{array}{ll}
u_{t}-\tri u=\mathbb{P} \dive \tau, \\
\tau_t=D(u),
\end{array}
\right.
\eeno
where $ \mathbb{P} $ is the Leray project operator. Then, one can see that both $u$ and $\mathbb{P} \dive \tau$ satisfy the following wave equation with damping term:
\[W_{tt}-\tri W_t-\f12\tri W=0.\]
Based on the above dissipative structure of $u$ and $\mathbb{P} \dive \tau$, the author in \cite{Zhu2018} proves
the global existence of the strong solution for the Oldroyd-B model with small initial data. However, from the linearized system, we can not obtain any dissipation for $\tau$. Thus the nonlinear term $(\tr\tau) \tau$ is a difficult term even though the initial data is small. In order to deal with this difficult term, we use the Lagrange coordinate to write the equation of $\tr \tau$, and then we can obtain some decay estimate of the quantity $\tr \tau$. In the later, we will see that the estimate of $\tr \tau$ is depend on some power of $\frac{1}{c_0}$ and the initial data. Since the initial data is small, it follows that $\frac{1}{c_0}$ is very large. The key point in this paper is to control the power of $\frac{1}{c_0}$ such that $\tr \tau$ is still small for any time. Combining the decay estimate of $\tr \tau$ with the dissipative structure for the linearized system of \eqref{PTT}, we can prove the global existence.

Throughout of the paper, we assume that $a=\lam=0$ and $b=\mu=\mu_1=\mu_2=1$. There is no derivative in the  additional term $(\tr\tau) \tau$ in \eqref{PTT}, the proof of local well-posedness for \eqref{PTT} is similar to the Oldroyd-B model(See \cite{Fernandez-Cara,Lin-Liu-Zhang2005}) and we omit the detail here.

Our main results can be stated as follow:

\begin{theo}\label{breaking}
Let $(u_0,\tau_0)\in H^2(\mathbb T^3)$, assume that $\tr \tau_0(x_0)<0$ for some $x_0\in \mathbb{T}^3$, then the corresponding solution to the system \eqref{PTT} blows up in finite time. More precisely, there exists a time $T$ with $0<T\leq -\f{1}{\tr \tau_0(x_0)}$ such that\beq
\lim_{t\to T^-} \tr\tau(t,q(t,x_0))=-\infty.
\eeq
\end{theo}

\begin{theo}\label{them}
  Suppose that $\dive u_0=0$, $\int_{\mathbb T^3}u_0dx=0$, $(\tau_0)_{ij}=(\tau_0)_{ji}$, and the initial data $(u_0,\tau_0)\in H^2(\mathbb T^3)$.  There exists $\ep_0,\tilde \ep_0$ such that if
 \beno
 0<\delta_0^2:=\|(u_0,\tau_0)\|^2_{ H^2(\mathbb T^3)}\leq \ep_0,
\eeno
\beno
 0<c_0^{-2}\|\na^2\tr\tau_0\|^2_{L^2(\mathbb T^3)}\leq \tilde\ep_0, \quad  \tr \tau_0(x)\geq c_0:=\f12\delta_0,
 \eeno
then the problem \eqref{PTT} admits a unique global solution $(u(t), \tau(t))$ satisfying that for all $t\geq0$: \beno
\begin{split}
\quad&\|u(t)\|_{H^2(\mathbb T^3)}^2+\|\tau(t)\|_{H^2(\mathbb T^3)}^2+\int_0^t\Big (\|\na u(s)\|_{H^2(\mathbb T^3)}^2+\|\mathbb P\dive\tau (s)\|_{H^1(\mathbb T^3)}^2+\|\na\tr\tau(s)\|_{H^1}^2 \Big)ds
\leq C\d_0^2,
\end{split}
\eeno
where $C>0$ is a positive constant independent of $t$.
Moreover, we have
\beq
\begin{split}
&\|u\|_{H^2(\mathbb T^3)}+\|\mathbb P\dive\tau\|_{L^2(\mathbb T^3)}\leq C(\|u_0\|_{H^2(\mathbb T^3)},\|\tau_0\|_{H^2(\mathbb T^3)})(1+t)^{-\f32+\f\ep2},
\end{split}
\eeq
where $\ep>0$ is a small constant.
\end{theo}

\begin{rema}
By virtue of Sobolev's embedding, one can see that
\[c_0\leq \tr \tau_0\leq 3\| \tau_0\|_{L^\infty}\leq 3\|\tau_0\|_{H^2}=3\delta_0,\]
it follows that $c_0=k\delta_0$ for some $k\leq 3$. For simplicity, we choose $k=\f12$ to prove our main result.
\end{rema}

\begin{rema}
The conditions $c_0^{-2}\|\na^2\tr\tau_0\|^2_{L^2(\mathbb T^3)}\leq \tilde\ep_0$ and $\tr \tau_0\geq c_0$ are feasible. For example, choose a periodic stress tensor $\tau_0$ such that
\[\tr \tau_0= \frac{3}{2}c_0+\frac{1}{4}c^2_0\tilde\ep_0\sin(x_1)\sin(x_2)\sin(x_3).\]
As $c_0>0$ is a small constant, one can easily check that
\beno
 0<c_0^{-2}\|\na^2\tr\tau_0\|^2_{L^2(\mathbb T^3)}\leq \tilde\ep_0, \quad  \tr \tau_0(x)\geq c_0.
 \eeno
\end{rema}

\begin{rema}
If $x\in \mathbb{R}^3$ and $\tau\in H^2(\mathbb{R}^3)$, we see that $\lim_{|x|\rightarrow \infty}\tau =0$. Thus, we can't suppose that $\tr \tau_0\geq c_0$ for some positive constant $c_0$. The Theorem \ref{them} is false in this case. The global wellposedness of \eqref{PTT} in the whole space is also an interesting problem.
\end{rema}

\begin{rema}
For the periodic Navier-Stokes equation, if $\int_{\mathbb T^3}u_0=0$, one can obtain the $H^2$-norm decay rate of velocity $u$ is $e^{-ct}$. However, for the system \eqref{PTT}, there is no any decay for the stress tensor $\tau$, and we only obtain an algebra decay in Theorem \ref{them}.
\end{rema}

The remainder of the paper is organized as follows. In Section 2 we give some key lemmas which will be used in the sequel. In Section 3, we prove the blow up phenomenon of the system \eqref{PTT}. Section 4 is devoted to study about the global existence and decay rate of the strong solution with small initial data.

\section{Preliminary}
We recall a few estimates for the flow of a smooth vector field.
\begin{lemm}\cite{Bahouri2011}\label{qt}
{\sl
Let $v$ be a smooth time-dependent vector field with bounded first order space derivatives. Let $q_t$ satisfy\beno
q_t(x)=x+\int_0^t v(s, q_s(x))ds.
\eeno
Then the flow $q_t$ is a $C^1$ diffeomorphism over $\mathbb T^3$, and for all $t\in\R^+$ we have \beno
\|\na q_t^{\pm 1}\|_{L^\infty}\leq \exp V(t),
\eeno
\beno
\|\na q_t^{\pm 1}-\text{Id}\|_{L^\infty}\leq \exp V(t)-1,
\eeno
\beno
\|\na^2 q_t^{\pm 1}\|_{L^\infty}\leq \exp V(t)\int_0^t \|\na^2 v(s)\|_{L^\infty}\exp V(s)ds,
\eeno
where \beno
V(t)=\int_0^t \|\na v(s)\|_{L^\infty}ds.
\eeno }
\end{lemm}


\begin{lemm}\label{lem}
For any smooth tensor $[\tau_{ij}]_{3\times3}$ and three dimensional vector $u$, it always holds that
\beno
\begin{split}
&\mathbb P\dive (u\cdot \na\tau)
=\mathbb P(u\cdot \na \mathbb P\dive\tau)+\mathbb P(\na u\cdot \na\tau)-\mathbb P(\na u\cdot \na \tri^{-1}\dive\dive\tau),\\
&\mathbb P\dive ((\tr\tau)\tau)
=\mathbb P((\tr\tau) \mathbb P\dive\tau)+\mathbb P(\tau\cdot\na (\tr\tau))-\mathbb P(\na (\tr\tau)\tri^{-1}\dive\dive\tau),
\end{split}
\eeno
where \beno
\begin{split}
&(\na u\cdot \na\tau)_i=\sum_j \pa_j u\cdot \na \tau_{ij},\quad
(\na u\cdot \na \tri^{-1}\dive\dive\tau)_i=\pa_i u\cdot \na \tri^{-1}\dive\dive\tau.
\end{split}
\eeno
\begin{proof}
The first equality has been proved in \cite{Zhu2018}. We only deal with the second equality. Using direct computation we have\beno
\mathbb P\dive ((\tr\tau)\tau)=\mathbb P(\tau\cdot\na (\tr\tau)+(\tr\tau)\dive \tau)=\mathbb P(\tau\cdot\na (\tr\tau))+\mathbb P((\tr\tau)\dive \tau).
\eeno
Denote $\mathbb P^\bot=\tri^{-1}\na\dive$, we compute as follows\beno
\begin{split}
\mathbb P((\tr\tau)\dive \tau)=&\mathbb P((\tr\tau)\mathbb P\dive \tau)+\mathbb P((\tr\tau)\tri^{-1}\na\dive\dive \tau)\\
=&\mathbb P((\tr\tau)\mathbb P\dive \tau)+\mathbb P\na((\tr\tau)\tri^{-1}\dive\dive \tau)-\mathbb P(\na(\tr\tau)\tri^{-1}\dive\dive \tau)\\
=&\mathbb P((\tr\tau)\mathbb P\dive \tau)-\mathbb P(\na(\tr\tau)\tri^{-1}\dive\dive \tau).
\end{split}
\eeno
Hence we proof the Lemma.
\end{proof}
\end{lemm}

By a directly calculation,  the following technical lemma holds true.
\begin{lemm}\label{lemm}
{\sl Let $r>0, c_0>0$, one has\beno
\begin{split}
\int_0^{\f t 2}e^{-(t-s)}(1+c_0s)^{-r} ds
\lesssim e^{-\f t2}\int_0^{\f t 2}(1+c_0s)^{-r} ds
\lesssim&{
\left\{\begin{array}{l}
c_0^{-1}e^{-\f t2}, \quad \text{for} \quad r>1,\\
c_0^{-1}e^{-\f t2}(1+c_0t)^{\ep},\quad \text{for} \quad r=1,\\
c_0^{-1}e^{-\f t2}(1+c_0t)^{1-r},\quad \text{for} \quad r<1,
\end{array}\right.}
\end{split}
\eeno
and\beno
\int_{\f t 2}^te^{-(t-s)}(1+c_0s)^{-r} ds
\lesssim(1+\f {c_0t} 2)^{-r}\int_{\f t 2}^te^{-(t-s)}ds
\lesssim(1+c_0t)^{-r},
\eeno
where $\epsilon>0$ is a small but fixed constant.
}
\end{lemm}

{\bf Notation.} Since all function spaces in through out the paper are over $\mathbb{T}^3$, for simplicity, we drop $\mathbb{T}^3$ in the notation of function spaces if there is no ambiguity.  $A\lesssim B$ stands for $A\leq C B$ for some constant $C>0$ independent of $A$ and $B$.

\section{Proof of the Theorem 1.1}
In this section, we are going to prove that the solution will blow up in finite time.
\begin{proof}
Suppose that $T^*$ is the maximal existence time for the solution $(u,\tau)$ with the initial data $u_0$ and $\tau_0$.  Applying operator $\tr$ to the second equation of \eqref{PTT}, and using the fact that $\tr D(u)=\dive u=0$ and $\tr Q(u,\tau)=0$,  we get \beq\label{tr}
\tr \tau_t+u\cdot \nabla \tr \tau+(\tr \tau)^2=0.
\eeq
Let us consider the trajectory equation,\beno
\f{d}{dt} q(t,x)=u(t, q(t,x)),\quad q(0,x)=x.
\eeno
Thanks to \eqref{tr}, we have\beno
\f{d}{dt} (\tr \tau)(t,q(t,x))=-(\tr \tau)^2(t,q(t,x)),
\eeno
which yields that\beq\label{tau}
\tr \tau(t,q(t,x))=\f{\tr \tau_0(x)}{1+\tr\tau_0(x)t}, \quad \forall t\in[0,T],
\eeq
where $0<T<T^*$. Since $\tr \tau_0(x_0)<0$, it follows that $\tr \tau(t,q(t,x_0))\rightarrow -\infty$ as $t\rightarrow -\frac{1}{\tr \tau_0(x_0)}$. Hence the maximal existence time $T^*\leq -\f{1}{\tr \tau_0(x_0)}$.
\end{proof}

\begin{rema}
Thanks to
$$\lim_{t\to T^{*-}} \tr\tau(t,q(t,x_0))=-\infty, $$
we have $\lim_{t\to T^{*-}} \|\tr \tau\|_{L^\infty}=\infty$. Using the Sobolev embedding, we see that
\[\lim_{t\to T^{*-}} \|\tau \|_{H^2}\geq \lim_{t\to T^{*-}} C\|\tau \|_{L^\infty} \geq \lim_{t\to T^{*-}} C\|\tr \tau \|_{L^\infty}=\infty.\]
Thus, our theorem implies that the $H^2$-norm of $\tau $ will blow up in finite time. However, we do not know whether the velocity will blow up or not before $T^*$.
\end{rema}

\begin{rema}
If $a\neq 0$ and $b>0$, by the same token we have
\beno
\f{d}{dt} (\tr \tau)(t,q(t,x))=-b(\tr \tau)^2(t,q(t,x))-a \tr \tau(t,q(t,x)) \leq -\f {b}{2}(\tr \tau)^2(t,q(t,x))+\f{a^2}{2 b}.
\eeno
From the above estimate, we can also obtain a blow up result when the initial data satisfy that $\tr \tau_0(x_0)<-\f{|a|}{b}$.
\end{rema}

\begin{rema}
Indeed we can obtain the blow up rate. Since
\beno
\f{d}{dt} (\tr \tau)(t,q(t,x))=-(\tr \tau)^2(t,q(t,x)),
\eeno
it follows that
\beno
\f{d}{dt} \f{1}{\tr \tau(t,q(t,x))}=-1,
\eeno
which yields that $\lim_{t\to T^{*-}}\tr\tau(t, q(t,x)))(T-t)=-1$.
\end{rema}

\section{Proof of the Theorem 1.2}
Now we turn our attention to prove the global existence and decay rate of the strong solution with small initial data and $\tr \tau_0(x)\geq c_0>0$.  First, we give the some basic energies and time-weighted energies as follows,
\beq
\mathcal E(0)=\|u_0\|_{H^2(\mathbb T^3)}^2+\|\tau_0\|_{H^2(\mathbb T^3)}^2,
\qquad\tilde{\mathcal E}(0)=c_0^{-2}\|\na^2\tr\tau_0\|_{L^2(\mathbb T^3)}^2,
\eeq
\beq
\mathcal E_1(t)=\sup_{0\leq s \leq t}(\|u(s)\|_{H^2(\mathbb T^3)}^2+\|\tau(s)\|_{H^2(\mathbb T^3)}^2)
+\int_0^t \|\na u(s)\|_{H^2(\mathbb T^3)}^2ds+\int_0^t \|\mathbb P \dive\tau(s)\|_{H^1(\mathbb T^3)}^2 ds,
\eeq
\beq
\begin{split}
\mathcal E_2(t)=&\sup_{0\leq s \leq t}(1+c_0s)^{3-\ep}(\|\na^2 u(s)\|_{L^2(\mathbb T^3)}^2+\|\na\mathbb P \dive\tau(s)\|_{L^2(\mathbb T^3)}^2)\\
&\quad+\int_0^t (1+c_0s)^{3-\ep}\|\na^3u(s)\|_{L^2(\mathbb T^3)}^2ds+\int_0^t (1+c_0s)^{3-\ep}\|\na\mathbb P \dive\tau(s)\|_{L^2(\mathbb T^3)}^2ds,
\end{split}
\eeq
\beq
\mathcal E_3(t)=\sup_{0\leq s \leq t}(1+s)^{3-\ep}(\|u(s)\|^2_{H^2(\mathbb T^3)}+\|\mathbb P \dive\tau(s)\|_{L^2(\mathbb T^3)}^2),
\eeq
\beq
\mathcal E_{4}(t)=c_0^{-1}\int_0^t (1+c_0s)^{3-\ep}\|\na\tr\tau(s)\|_{L^2(\mathbb T^3)}^2 ds,
\eeq
\beq
\mathcal E_{5}(t)=c_0^{-1}\int_0^t (1+c_0s)^{3-\ep}\|\na^2\tr\tau(s)\|_{L^2(\mathbb T^3)}^2 ds,
\eeq
for any $\ep>0$ is a fixed but small constant, and $\mathbb P=\mathbb I-\tri^{-1}\na\dive$ is the Leray projection operator.

We shall derive the a priori estimates of $\mathcal E_1(t)$, $\mathcal E_2(t)$, $\mathcal E_3(t)$, $\mathcal E_4(t)$ and $\mathcal E_5(t)$ respectively.

\subsection{The estimates of $\mathcal E_1(t)$}
Applying the operator $\na^{k}$ $(k=0,1,2)$ to the \eqref{PTT} system and by virtue of the standard energy estimate, we have\beq\label{e11}
\begin{split}
&\f12\f{d}{dt}\Big\{\|u\|_{H^2}^2+\|\tau\|_{H^2}^2\Big\}
+\|\na u\|_{H^2}^2\\
=&-\sum_{k=0}^2\int_{\mathbb T^3}\na^{k}(u\cdot \na u)\na^{k} udx-\sum_{k=0}^2\int_{\mathbb T^3}\na^k (u\cdot \na\tau+ (\tr\tau)\tau+Q(\tau,\na u))\na^k\tau dx\\
&\quad+\sum_{k=0}^2\int_{\mathbb T^3}(\na^{k}\dive \tau\na^{k} u+\na^k D(u)\na^k\tau)dx=I_1+I_2+I_3.
\end{split}
\eeq
For the first term $I_1$, notice that $\dive u=0$, by virtue of H\"{o}lder's inequality and Sobolev's embedding, we obtain that
\beno
\begin{split}
\int_0^t I_1(s)ds=&-\sum_{k=0}^2\int_0^t \int_{\mathbb T^3}\na^{k}(u\cdot \na u)\na^{k} u(s)dxds\\
\lesssim&\int_0^t \Big(\|\na u(s)\|_{L^2}\|\na u(s)\|_{L^3}\|\na u(s)\|_{L^6}+\|\na u(s)\|_{L^6}\|\na^2 u(s)\|_{L^3}\|\na^2 u(s)\|_{L^2}\Big)ds\\
\lesssim&\sup_{0\leq s \leq t}\|u(s)\|_{H^2}\int_0^t \|\na u(s)\|_{H^2}^2ds \lesssim \mathcal E_1^{\f32}(t).
\end{split}
\eeno
Because of $\int_{\mathbb T^3} u_0(x)dx=0$, we deduce from the first equation of \eqref{PTT} that $\int_{\mathbb T^3} u(t, x)dx= 0$. For the term $I_2$, by virtue of Poincar\'{e}'s inequality, we get that\beno
\begin{split}
\int_0^t I_2(s)ds
=&-\sum_{k=0}^2\int_0^t\int_{\mathbb T^3}\na^k (u\cdot \na\tau+ (\tr\tau)\tau+Q(\tau,\na u))\na^k\tau (s)dxds\\
\lesssim&\int_0^t \Big\{\big[\|\na u(s)\|_{L^\infty}\|\na \tau(s)\|_{H^1}^2+\|\na^2 u(s)\|_{L^3}\|\na \tau(s)\|_{L^6}\|\na^2 \tau(s)\|_{L^2}\big]\\
&\quad+ \big[(\|\tr\tau(s)\|_{L^\infty}+\|\na^2\tr\tau(s)\|_{L^2})\|\tau(s)\|_{H^2}^2+ \|\na\tr\tau(s)\|_{L^3}\|\tau(s)\|_{W^{1,6}}\|\na \tau(s)\|_{H^1}\big]\\
&\quad+\big[(\|\na u(s)\|_{L^\infty}+\|\na^3 u(s)\|_{L^2})\|\tau(s)\|_{H^2}^2+ \|\na^2 u(s)\|_{L^3}\|\tau(s)\|_{W^{1,6}}\|\na \tau(s)\|_{H^1}\big]\Big\}ds\\
\lesssim &\sup_{0\leq s \leq t}\|\tau(s)\|_{H^2}^2\int_0^t \Big\{\|\na^3 u(s)\|_{L^2}+\|\na\tr\tau(s)\|_{L^2}^\f12\|\na^2\tr\tau(s)\|_{L^2}^\f12+\|\na^2\tr\tau(s)\|_{L^2}\Big\}ds\\
\lesssim &\sup_{0\leq s \leq t}\|\tau(s)\|_{H^2}^2\Big(\int_0^t (1+c_0s)^{-3+\ep}ds\Big)^{\f12}\\
&\quad\times\Big(\int_0^t (1+c_0s)^{3-\ep}(\|\na^3 u(s)\|_{L^2}^2+\|\na\tr\tau(s)\|_{L^2}\|\na^2\tr\tau(s)\|_{L^2}+\|\na^2\tr\tau(s)\|_{L^2}^2)ds\Big)^{\f12}\\
\lesssim &c_0^{-\f12}\mathcal E_1(t)\mathcal E_2^{\f12}(t)+\mathcal E_1(t)\mathcal E_{4}^{\f14}(t)\mathcal E_{5}^{\f14}(t)+\mathcal E_1(t)\mathcal E_{5}^{\f12}(t).
\end{split}
\eeno
For the last term $I_3$, using integration by parts and the symmetry $\tau_{ij}=\tau_{ji}$, we have\beno
I_3=\sum_{k=0}^2\int_{\mathbb T^3}(\na^{k}\dive \tau\na^{k} u+\na^k D(u)\na^k\tau)dx=0.
\eeno
Integrating \eqref{e11} with time, according to above estimates, we deduce that
\beq\label{E11}
\begin{split}
\mathcal E_{11}=&\sup_{0\leq s \leq t}(\|u(s)\|_{H^2}^2+\|\tau(s)\|_{H^2}^2)
+\int_0^t \|\na u(s)\|_{H^2}^2ds\\
\lesssim& \mathcal E_0(t)+c_0^{-\f12}\Big(\mathcal E_1^{\f32}(t)+\mathcal E_2^{\f32}(t)\Big)+\mathcal E_1(t)\mathcal E_{4}^{\f14}(t)\mathcal E_{5}^{\f14}(t)+\mathcal E_1(t)\mathcal E_{5}^{\f12}(t).
\end{split}
\eeq

Operating $\mathbb P$ on the first equation of the \eqref{PTT} system, we have the following equation\beno
u_{t}+\mathbb P(u\cdot \na u)-\tri u=\mathbb P\dive \tau.
\eeno
Applying $\na^k$ $(k=0,1)$ to the above equation,  taking inner product with $\na^{k}\mathbb P\dive \tau $, then we obtain that \beq\label{e12}
\begin{split}
\|\mathbb P\dive \tau\|_{H^1}^2=&\sum_{k=0,1}\int_{\mathbb T^3}\na^ku_{t}\na^k\mathbb P\dive \tau dx+\sum_{k=0,1}\int_{\mathbb T^3}\na^k\mathbb P(u\cdot \na u)\na^k\mathbb P\dive \tau dx\\
&\quad-\sum_{k=0,1}\int_{\mathbb T^3}\na^k\tri u\na^k\mathbb P\dive \tau dx=J_1+J_2+J_3.
\end{split}
\eeq
For the first term $J_1$, using integration by parts, we rewrite this term into the following form\beno
\begin{split}
J_1=&\sum_{k=0,1}\f{d}{dt}\Big\{\int_{\mathbb T^3} \na^ku\na^k\mathbb P\dive \tau dx\Big\}-\sum_{k=0,1}\int_{\mathbb T^3} \na^ku\na^k\mathbb P\dive \tau_t dx=J_{11}+J_{12},
\end{split}
\eeno
where \beno
\int_0^t J_{11} (s)ds\lesssim \int_0^t \f{d}{ds}\Big\{\|u(s)\|_{H^1}\|\tau (s)\|_{H^2}\Big\} ds
\lesssim \mathcal E(0)+\mathcal E_{11}(t).
\eeno
Operating $\mathbb P\dive$ on the second equation of the \eqref{PTT} system, we have the following equation\beno
\mathbb P\dive\tau_t+\mathbb P\dive(u\cdot \na\tau+ (\tr\tau)\tau+Q(\tau,\na u))=\f12\tri u,
\eeno
which leads to
\beno
\begin{split}
\int_0^t J_{12}(s)ds
=&\sum_{k=0,1}\int_0^t\int_{\mathbb T^3} (-1)^{k}\tri^{k}u(\mathbb P\dive(u\cdot \na\tau+ (\tr\tau)\tau+Q(\tau,\na u))-\f12\tri u) (s)dxds\\
\lesssim&\int_0^t \Big\{\big[\|\na u\|^2_{L^2}+\|u(s)\|_{L^\infty}\|u(s)\|_{L^2}\|\na \mathbb P\dive\tau(s)\|_{L^2}+\|u(s)\|_{L^\infty}\|\na u(s)\|_{L^2}\|\tau(s)\|_{H^1}\\
&\quad+\|u(s)\|_{L^\infty}\|\tr \tau(s)\|_{L^2}\|\mathbb P\dive\tau(s)\|_{L^2}+\|u(s)\|_{L^\infty}\|\na \tr\tau(s)\|_{L^2}\|\tau(s)\|_{L^2}\big]\\
&\quad+\big[\|\na^2 u\|^2_{L^2}+\|\na^2 u\|_{L^2}\|u(s)\|_{L^\infty}\|\na \mathbb P\dive\tau(s)\|_{L^2}+\|\na^2 u\|_{L^2}\|\na u(s)\|_{L^\infty}\|\na \tau(s)\|_{L^2}\\
&\quad+\|\na^2 u(s)\|_{L^2}\|\tr \tau(s)\|_{L^\infty}\|\mathbb P\dive\tau(s)\|_{L^2}\\
&\quad+\|\na^2u(s)\|_{L^2}(\|\na \tr\tau(s)\|_{L^2}+\|\na u(s)\|_{L^2})\|\tau(s)\|_{L^\infty}
\big]\Big\}ds\\
\lesssim&\int_0^t \|\na u(s)\|_{H^1}^2 ds+\sup_{0\leq s\leq t}(\|u(s)\|_{H^2}+\|\tau(s)\|_{H^2})\\
&\quad\times\int_0^t (\|\na u(s)\|_{H^2}^2+\|\mathbb P\dive\tau(s)\|_{H^1}^2+\|\na u(s)\|_{H^2}\|\na \tr\tau(s)\|_{L^2})ds\\
\lesssim& \mathcal E_{11}+\mathcal E_1^{\f32}(t)+c_0^\f12\mathcal E_1(t)\mathcal E_4^{\f12}(t).
\end{split}
\eeno
Thus, we get that\beno
\int_0^t J_1 (s)ds\lesssim \mathcal E(0)+\mathcal E_{11}+\mathcal E_1^{\f32}(t)+c_0^\f12\mathcal E_1(t)\mathcal E_4^{\f12}(t).
\eeno
For the second term $J_2$ and the last term $J_3$, we can directly compute\beno
\begin{split}
\int_0^t J_2 (s)ds=& \sum_{k=0,1}\int_0^t \int_{\mathbb T^3}\na^k\mathbb P(u\cdot \na u)\na^k\mathbb P\dive \tau (s)dxds\\
\lesssim& \int_0^t (\|u(s)\|_{L^\infty}\|\na u(s)\|_{H^1}+\|\na u(s)\|_{L^\infty}\|\na u(s)\|_{L^2})\| \mathbb P\dive \tau(s)\|_{H^1}ds\\
\lesssim&\sup_{0\leq s\leq t}\|u(s)\|_{H^2}\Big(\int_0^t \|\na u(s)\|^2_{H^2}ds\Big)^{\f12}\Big(\int_0^t \|\mathbb P\dive\tau(s)\|^2_{H^1}ds\Big)^{\f12}
\lesssim \mathcal E_1^{\f32}(t),
\end{split}
\eeno
and\beno
\begin{split}
\int_0^t J_3 (s)ds= -\sum_{k=0,1}\int_0^t \int_{\mathbb T^3}\na^k\tri u\na^k\mathbb P\dive \tau (s)dxds
\lesssim \mathcal E_{11}^\f12(t)\mathcal E_{12}^\f12(t).
\end{split}
\eeno
Integrating \eqref{e12} with time, combining the above estimates, yields that\beq\label{E12}
\mathcal E_{12}(t)=\int_0^t \|\mathbb P \dive\tau(s)\|_{H^1}^2 ds\lesssim
 \mathcal E(0)+\mathcal E_{11}+\mathcal E_1^{\f32}(t)+c_0^\f12\mathcal E_1(t)\mathcal E_4^{\f12}(t).
\eeq

Adding up $C_1\times$\eqref{E11} and \eqref{E12}, choosing $C_1>0$ fixed but large enough, then we obtain that\beq\label{E1}
\begin{split}
\mathcal E_1(t)\lesssim& \mathcal E(0)+c_0^\f12\mathcal E_1(t)\mathcal E_4^\f12(t)+ c_0^{-\f12}\Big(\mathcal E_1^{\f32}(t)+\mathcal E_2^{\f32}(t)\Big)+\mathcal E_1(t)\mathcal E_{4}^{\f14}(t)\mathcal E_{5}^{\f14}(t)+\mathcal E_1(t)\mathcal E_{5}^{\f12}(t).
\end{split}
\eeq

\subsection{The estimates of $\mathcal E_2(t)$}
Taking $\na^2$ on the first equation of the system \eqref{PTT} and applying $\na \mathbb P\dive$ on the second equation of the system \eqref{PTT}, then we obtain the following system:
\begin{align}\label{OB1}
\left\{
\begin{array}{ll}
\na^{2}u_{t}+\na^{2}(u\cdot \na u)-\na^{2}\tri u+\na^{2}\na p=\na^{2}\dive \tau, \\[1ex]
\na \mathbb P\dive\tau_t+\na \mathbb P\dive(u\cdot \na\tau+ (\tr\tau)\tau+Q(\tau,\na u))=\f12\na \tri u. \\[1ex]
\end{array}
\right.
\end{align}
Taking inner product with $\na^{2} u$ for the first equation of the system \eqref{OB1}, taking inner product with $2\na \mathbb P\dive\tau$ for the second equation of the system \eqref{OB1}, and adding the time weight $(1+c_0t)^{3-\ep}$, then we have\beq\label{e21}
\begin{split}
&\f12\f{d}{dt}\Big\{(1+c_0t)^{3-\ep}(\|\na^2 u\|_{L^2}^2+2\|\na \mathbb P\dive\tau\|_{L^2}^2)\Big\}
+(1+c_0t)^{3-\ep}\|\na^3 u\|_{L^2}^{2}\\
=&\f{{(3-\ep)c_0}}2(1+c_0t)^{2-\ep}(\|\na^2 u\|_{L^2}^2+2\| \na\mathbb P\dive\tau\|_{L^2}^2)\\
&\quad-(1+c_0t)^{3-\ep}\int_{\mathbb T^3}\na^{2}(u\cdot \na u)\na^{2} udx\\
&\quad-2(1+c_0t)^{3-\ep}\int_{\mathbb T^3}\na \mathbb P\dive(u\cdot \na\tau+ (\tr\tau)\tau+Q(\tau,\na u))\na \mathbb P\dive\tau dx\\
&\quad+(1+c_0t)^{3-\ep}\int_{\mathbb T^3}(\na^{2}\dive \tau\na^{2} u+\na \tri u\na\mathbb P\dive\tau)dx=K_1+K_2+K_3+K_4.
\end{split}
\eeq
For the first term $K_1$ and the second term $K_2$, by virtue of Poincare's inequality, we can directly derive\beno
\begin{split}
\int_0^t K_1(s) ds
=&\int_0^t \f{{(3-\ep)c_0}}2(1+c_0s)^{2-\ep}(\|\na^2 u(s)\|_{L^2}^2+2\|\na \mathbb P\dive\tau(s)\|_{L^2}^2) ds\\
\lesssim&c_0\int_0^t (1+c_0s)^{3-\ep}(\|\na^3 u(s)\|_{L^2}^2+\|\na \mathbb P\dive\tau(s)\|_{L^2}^2) ds
\lesssim c_0\mathcal E_2(t),
\end{split}
\eeno
and\beno
\begin{split}
\int_0^t K_2(s) ds=&-\int_0^t \int_{\mathbb T^3}(1+c_0s)^{3-\ep}\na^2(u\cdot \na u)\na^2u(s)dx ds\\
\lesssim&\int_0^t (1+c_0s)^{3-\ep}\|u(s)\|_{H^2}\|\na^3 u(s)\|_{L^2}^2ds
\lesssim \mathcal E_1^{\f12}(t) \mathcal E_2(t).
\end{split}
\eeno
We turn to deal with the wildest term $K_3$. Applying Lemma \ref{lem}, H\"{o}lder's inequality and Sobolev's embedding, then we deduce that\beno
\begin{split}
\int_0^t K_3(s) ds
=&-2\int_0^t\int_{\mathbb T^3}(1+c_0s)^{3-\ep}\na \mathbb P\dive(u\cdot \na\tau+ (\tr\tau)\tau+Q(\tau,\na u))\na \mathbb P\dive\tau (s)dxds\\
\lesssim&\int_0^t (1+c_0s)^{3-\ep}\Big[(\|\na u(s)\|_{L^\infty}+\|\tr\tau(s)\|_{L^\infty})\|\na^2 \tau(s)\|_{L^2}\|\na\mathbb P\dive\tau(s)\|_{L^2}\\
&\quad+(\|\na^2 u(s)\|_{L^3}+\|\na\tr\tau(s)\|_{L^3})\|\na \tau(s)\|_{L^6}\|\na\mathbb P\dive\tau(s)\|_{L^2}\\
&\quad+\|\tau(s)\|_{L^\infty}(\|\na^2 \tr \tau(s)\|_{L^2}+\|\na^3u(s)\|_{L^2})\|\na\mathbb P\dive\tau(s)\|_{L^2}\Big]ds\\
\lesssim&\sup_{0\leq s\leq t}\|\tau(s)\|_{H^2}\int_0^t (1+c_0s)^{3-\ep}(\|\na^3 u(s)\|_{L^2}\|\na\mathbb P\dive\tau(s)\|_{L^2}+\|\na \tr \tau(s)\|_{H^1}\|\na\mathbb P\dive\tau(s)\|_{L^2})ds\\
\lesssim& \mathcal E_1^{\f12}(t)\mathcal E_2(t) +c_0^\f12\mathcal E_1^{\f12}(t)\mathcal E_2^{\f12}(t)\Big(\mathcal E_4^{\f12}(t)+\mathcal E_5^{\f12}(t)\Big).
\end{split}
\eeno
Taking advantage of integration by parts and $\dive u=0$, we can compute\beno
K_4=(1+c_0s)^{3-\ep}\int_{\mathbb T^3}(\na^{2}\dive \tau\na^{2} u+\na\tri u\na\mathbb P\dive\tau)dx=0.
\eeno
Integrating \eqref{e21} with time, combing the above estimates, yields that\beq\label{E21}
\begin{split}
\mathcal E_{21}(t)=&\sup_{0\leq s \leq t}(1+c_0s)^{3-\ep}(\|\na^2 u(s)\|_{L^2}^2+\|\na\mathbb P \dive\tau(s)\|_{L^2}^2)+\int_0^t (1+c_0s)^{3-\ep}\|\na^3u(s)\|_{L^2}^2ds\\
\lesssim & \mathcal E(0)+c_0\mathcal E_2(t)+\mathcal E_1^{\f32}(t)+\mathcal E_2^{\f32}(t) +c_0^\f12\mathcal E_1^{\f12}(t)\mathcal E_2^{\f12}(t)\Big(\mathcal E_4^{\f12}(t)+\mathcal E_5^{\f12}(t)\Big).
\end{split}
\eeq

Operating $\na \mathbb P$ on the first equation of the \eqref{PTT} system, we have the following equation\beno
\na u_{t}+\na \mathbb P(u\cdot \na u)-\na \tri u=\na \mathbb P\dive \tau.
\eeno
Taking $L^2$ inner product with $\na \mathbb P\dive \tau$, adding the time weight $(1+c_0t)^{3-\ep}$, then we obtain that \beq\label{e22}
\begin{split}
&(1+c_0t)^{3-\ep}\|\na \mathbb P\dive \tau\|_{L^2}\\
=&(1+c_0t)^{3-\ep}\int_{\mathbb T^3} \na u_{t}\na \mathbb P\dive \tau dx+(1+c_0t)^{3-\ep}\int_{\mathbb T^3}\na \mathbb P(u\cdot \na u)\na \mathbb P\dive \tau dx\\
&\quad-(1+c_0t)^{3-\ep}\int_{\mathbb T^3}\na \tri u\na \mathbb P\dive \tau dx=L_1+L_2+L_3.
\end{split}
\eeq
For the first term $L_1$, integrating by parts, we see that \beno
\begin{split}
L_1=&\f{d}{dt}\Big\{(1+c_0t)^{3-\ep}\int_{\mathbb T^3} \na u\na \mathbb P\dive \tau dx\Big\}-(3-\ep)c_0(1+c_0t)^{2-\ep}\int_{\mathbb T^3} \na u\na \mathbb P\dive \tau dx\\
&\quad-(1+c_0t)^{3-\ep}\int_{\mathbb T^3} \na u\na \mathbb P\dive \tau_t dx=L_{11}+L_{12}+L_{13},
\end{split}
\eeno
where \beno
\int_0^t L_{11} (s)ds\leq \int_0^t \f{d}{ds}\Big\{(1+c_0s)^{3-\ep}\|\na u(s)\|_{L^2}\|\na \mathbb P\dive \tau (s)\|_{L^2}\Big\} ds
\leq \mathcal E(0)+\mathcal E_{21}(t),
\eeno
and\beno
\begin{split}
\int_0^t L_{12}(s)ds
\lesssim&c_0\Big(\int_0^t (1+c_0s)^{3-\ep}\|\na^3 u(s)\|^2_{L^2}ds\Big)^{\f12}
\Big(\int_0^t (1+c_0s)^{3-\ep}\|\na \mathbb P\dive \tau(s)\|^2_{L^2}ds\Big)^{\f12}
\lesssim c_0\mathcal E_{21}^\f12(t)\mathcal E_{22}^\f12(t).
\end{split}
\eeno
Operating $\mathbb P\dive$ on the second equation of the system \eqref{PTT}, we have the following equation\beno
\mathbb P\dive\tau_t+\mathbb P\dive(u\cdot \na\tau+ (\tr\tau)\tau+Q(\tau,\na u))=\f12\tri u,
\eeno
from which we deduce that\beno
\begin{split}
\int_0^t L_{13}(s)ds
=&-\int_0^t\int_{\mathbb T^3} (1+c_0s)^{3-\ep}\tri u (\mathbb P\dive(u\cdot \na\tau+ (\tr\tau)\tau+Q(\tau,\na u))-\f12\tri u) dxds.
\end{split}
\eeno
Applying the Lemma \ref{lem}, we get that\beno
\begin{split}
\int_0^t L_{13}(s)ds
\lesssim&\int_0^t (1+c_0s)^{3-\ep}\Big\{\|\na^2 u(s)\|_{L^2}^2+\|\na^2 u(s)\|_{L^6}\big[\|u(s)\|_{L^3}\|\na\mathbb P\dive\tau(s)\|_{L^2}+\|\na u(s)\|_{L^2}\|\na\tau(s)\|_{L^3}\\
&\quad+\|\tau(s)\|_{L^3}(\|\na \tr\tau(s)\|_{L^2}+\|\na^2 u(s)\|_{L^2})\big]ds+\|\na^2 u(s)\|_{L^2} \|\tr\tau(s)\|_{L^\infty}\|\mathbb P\dive\tau(s)\|_{L^2}\Big\},\\
\lesssim &\mathcal E_{21}(t)+\mathcal E_1^{\f12}(t)\mathcal E_2(t)
+c_0^\f12\mathcal E_1^{\f12}(t)\mathcal E_2^{\f12}(t)\Big(\mathcal E_4^{\f12}(t)+\mathcal E_5^{\f12}(t)\Big).
\end{split}
\eeno
From the above estimates, we obtain that\beno
\int_0^t L_1 (s)ds\lesssim \mathcal E(0)+\mathcal E_{21}(t)+c_0\mathcal E_{21}^\f12(t)\mathcal E_{22}^\f12(t)+\mathcal E_1^{\f32}(t)+\mathcal E_2^{\f32}(t)
+c_0^\f12\mathcal E_1^{\f12}(t)\mathcal E_2^{\f12}(t)\Big(\mathcal E_4^{\f12}(t)+\mathcal E_5^{\f12}(t)\Big).
\eeno
For the second term $L_2$ and the last term $L_3$, we can directly derive\beno
\begin{split}
\int_0^t L_2 (s)ds
\lesssim& \int_0^t (1+c_0s)^{3-\ep}(\|u\|_{L^\infty}\|\na^2 u\|_{L^2}+\|\na u\|_{L^\infty}\|\na u\|_{L^2})\|\na \mathbb P\dive \tau\|_{L^2}ds
\lesssim \mathcal E_1^{\f12}(t)\mathcal E_2(t),
\end{split}
\eeno
and\beno
\begin{split}
\int_0^t L_3 (s)ds
\lesssim& \int_0^t (1+c_0s)^{3-\ep}\|\na^3 u\|_{L^2}\|\na \mathbb P\dive \tau\|_{L^2}ds
\lesssim \mathcal E_{21}^\f12(t)\mathcal E_{22}^\f12(t).
\end{split}
\eeno
Integrating \eqref{e22} with time, combining the above estimates, then we have\beq\label{E22}
\begin{split}
\mathcal E_{22}(t)=&\int_0^t (1+c_0s)^{3-\ep}\|\na\mathbb P \dive\tau(s)\|_{L^2}^2ds\\
\lesssim & \mathcal E(0)+\mathcal E_{21}(t)+\mathcal E_1^{\f32}(t)+\mathcal E_2^{\f32}(t) +c_0^\f12\mathcal E_1^{\f12}(t)\mathcal E_2^{\f12}(t)\Big(\mathcal E_4^{\f12}(t)+\mathcal E_5^{\f12}(t)\Big).
\end{split}
\eeq

Adding up $C_2\times$\eqref{E21} and  \eqref{E22}, choosing $C_2>0$ is fixed but large enough, then we obtain that\beq\label{E2}
\mathcal E_2(t)\lesssim \mathcal E(0)+c_0^\f12\mathcal E_1^{\f12}(t)\mathcal E_2^{\f12}(t)\Big(\mathcal E_4^{\f12}(t)+\mathcal E_5^{\f12}(t)\Big)+\mathcal E_1^{\f32}(t)+\mathcal E_2^{\f32}(t).
\eeq

\subsection{The estimates of $\mathcal E_3(t)$}
The main difficult is to estimate $\mathcal E_3(t)$.
Let us consider the initial value problem for the following system:
\begin{align}\label{ob}
\left\{
\begin{array}{ll}
u_{t}-\tri u-\mathbb P\dive \tau=-\mathbb P(u\cdot \na u), \quad (t, x) \in \mathbb{R}^{+}\times \mathbb{T}^3, \\[1ex]
(\mathbb P\dive \tau)_t-\f12\tri u=-\mathbb P\dive (u\cdot \na\tau+ (\tr\tau)\tau+Q(\tau,\na u)),  \quad (t, x) \in \mathbb{R}^{+}\times \mathbb{T}^3, \\[1ex]
\dive u=0,  \quad (t, x) \in \mathbb{R}^{+}\times \mathbb{T}^3, \\[1ex]
u|_{t=0}=u_{0}(x),  \quad \tau|_{t=0}=\tau_0(x), \quad x \in  \mathbb{T}^3.\\[1ex]
\end{array}
\right.
\end{align}
In terms of the semigroup theory, the solution $(u, \mathbb P\dive \tau)$
of the \eqref{ob} can be expressed for $U = (u, \mathbb P\dive \tau)^ t $ as \beno
U_t(t,x)= BU(t,x)+H(t,x), \quad U(0,x) = U_0(x), \quad t \geq 0,
\eeno
where $H(t,x)=(H_1(t,x),H_2(t,x))^t=(-\mathbb P(u\cdot \na u)(t,x), -\mathbb P\dive (u\cdot \na\tau+ (\tr\tau)\tau+Q(\tau,\na u))(t,x))^t$.
Applying the Fourier transform to the system \eqref{ob},
we have \beno
\widehat{U}_t(t,k) = A(k)\widehat{U}(t,k)+\widehat{H}(t,k), \quad \widehat{U}(0,k) = \widehat{U}_0(k), \quad t \geq 0,
\eeno
where $k=(k_1,k_2,k_3)\in\mathbb Z^3$, $A(k)$ is defined as \beno
A(k) = {\left(
\begin{matrix}
  -|k|^2I_{3\times3} &  I_{3\times3}\\
 -\f{|k|^2}2I_{3\times3} & 0
\end{matrix}
\right).}
\eeno
The eigenvalues of the matrix $A(k)$ are computed from the determinant \beno
\det(A(k)-\lambda I) = \Big(\lambda^2+ |k|^2\lambda + \f{|k|^2}2\Big)^3=0,
\eeno
which implies \beno
\lambda_{1}=\lambda_1(|k|)(\text{triple root}), \quad \lambda_{2}=\lambda_2(|k|)(\text{triple root}).
\eeno
Thus, the semigroup $e^{tA(k)}$ is expressed as \beno
e^{tA(k)}=e^{\lambda_{1} t}P_{1}+e^{\lambda_{2} t}P_{2},
\eeno
where the project operators $P_i$ can be computed as \beno
P_i=\prod_{i\ne j}\f{A(k)-\lambda_j I}{\lambda_i-\lambda_j}.
\eeno
By a direct computation, we deduce the Fourier transform $\widehat G(t,k)$ of
Green's function $G(t,x)$ of $e^{tB}$ as \beno
\begin{split}
\widehat G(t,k)&=e^{tA(k)}
={\left(
\begin{matrix}
 (\f{\lambda_1e^{\lambda_2 t}-\lambda_2e^{\lambda_1 t}}{\lambda_1-\lambda_2}-\f{|k|^2(e^{\lambda_1 t}-e^{\lambda_2 t})}{\lambda_1-\lambda_2})I_{3\times3} & \f{e^{\lambda_1 t}-e^{\lambda_2 t}}{\lambda_1-\lambda_2}I_{3\times3} \\
 -\f{|k|^2(e^{\lambda_1 t}-e^{\lambda_2 t})}{2(\lambda_1-\lambda_2)}I_{3\times3} & \f{\lambda_1e^{\lambda_2 t}-\lambda_2e^{\lambda_1 t}}{\lambda_1-\lambda_2}I_{3\times3}
\end{matrix}
\right)}
={\left(
\begin{matrix}
\widehat N (t,k) \\
\widehat M(t,k)
\end{matrix}
\right).}
\end{split}
\eeno
Then we have the following decomposition for the Fourier transform of $(u, \mathbb P\dive \tau)$ as \beq
\begin{split}
&\widehat {u}(t,k) =  \widehat N(t,k)\cdot \widehat{U}_0 (k)+\int_0^t \widehat N (t-s,k)\cdot   \widehat H(s,k)ds,\\
&\widehat {\mathbb P\dive \tau} (t,k)=  \widehat M (t,k) \cdot \widehat{U}_0(k) +\int_0^t \widehat M (t-s,k)\cdot   \widehat H(s,k)ds.
\end{split}
\eeq

We need to verify the approximation of the Fourier transform of Green's function for all $k\in\mathbb Z^3$.\\
If $|k|\leq 1$, we see that \beq
\begin{split}
\lambda_1 = -\f12+ \f i 2,\qquad
\lambda_2 = -\f12 - \f i 2,
\end{split}
\eeq
which leads to \beno
\begin{split}
\f{\lambda_1e^{\lambda_2 t}-\lambda_2e^{\lambda_1 t}}{\lambda_1-\lambda_2}
=e^{-\f t2}\Big(\cos(\f t2)+\sin(\f t2)\Big)\lesssim e^{-\f t2},\qquad
\f{e^{\lambda_1 t}-e^{\lambda_2 t}}{\lambda_1-\lambda_2}
=2e^{-\f t2}\sin(\f t2)\lesssim e^{-\f t2}.
\end{split}
\eeno
If $|k|\geq 2$, we see that \beq
\begin{split}
\lambda_1 = -\f12 |k|^2 - \f12 \sqrt{(|k|^2-1)^2-1},\qquad
\lambda_2 = -\f12 |k|^2 + \f12 \sqrt{(|k|^2-1)^2-1},
\end{split}
\eeq
and we have \beno
\begin{split}
\f{\lambda_1e^{\lambda_2 t}-\lambda_2e^{\lambda_1 t}}{\lambda_1-\lambda_2}
\lesssim e^{-\f t2},\qquad
\f{e^{\lambda_1 t}-e^{\lambda_2 t}}{\lambda_1-\lambda_2}
\lesssim\f1{|k|^2}e^{-\f t2}.
\end{split}
\eeno

An easy computation together with the formula of the $\widehat G(t,k)$ gives \beno
\begin{split}
\widehat N (t,k)\cdot \widehat{U}_0(k)&=\Big( \f{\lambda_1e^{\lambda_2 t}-\lambda_2e^{\lambda_1 t}}{\lambda_1-\lambda_2}-\f{|k|^2(e^{\lambda_1 t}-e^{\lambda_2 t})}{\lambda_1-\lambda_2}\Big)\widehat{u}_0(k)+\f{e^{\lambda_1 t}-e^{\lambda_2 t}}{\lambda_1-\lambda_2}\widehat{\mathbb P\dive \tau}_0(k),
\end{split}
\eeno
and\beno
\begin{split}
\widehat M (t,k)\cdot \widehat{U}_0(k)&=-\f{|k|^2(e^{\lambda_1 t}-e^{\lambda_2 t})}{2(\lambda_1-\lambda_2)}\widehat{u}_0+\f{\lambda_1e^{\lambda_2 t}-\lambda_2e^{\lambda_1 t}}{\lambda_1-\lambda_2}\widehat{\mathbb P\dive \tau}_0.
\end{split}
\eeno
Since $\int_{\mathbb T^3} u(t, x)dx= 0$, $\int_{\mathbb T^3} \mathbb P\dive\tau(t,x)dx=0$, $\forall t\geq0$ it follows that $\widehat{U}(0)=0$. Then we have the $L^2$-decay rate on the derivatives of $(N \ast{U}_0, M \ast{U}_0)$ as \beq\label{e6-1}
\begin{split}
\| {\na^{\al}(N \ast{U}_0)}(t)\|^2_{L^2(\mathbb T^3)}
&\lesssim \sum_{k\in\mathbb Z^3}|\widehat N (t,k)\cdot \widehat{\na^{\al}{U}_0}(k)|^2\\
&\lesssim \sum_{k\in\mathbb Z^3/\{0\}}e^{-t}|\widehat{\na^{\al}{u}_0}(k)|^2+\sum_{k\in\mathbb Z^3/\{0\}}\f{e^{-t}}{|k|^2}|\widehat{\na^{\al}{{\tau}}_0}(k)|^2\\
&\lesssim e^{-t}\big(\|\na^{\al}{u}_0\|^2_{L^2(\mathbb T^3)} +\|\na^{\al}{\tau}_0\|^2_{L^2(\mathbb T^3)}\big),\quad \forall \al \in\N^3,
\end{split}
\eeq
and
\beq\label{e6-2}
\begin{split}
\| {\na^{\al}(M \ast{U}_0)}(t)\|_{L^2(\mathbb T^3)}^2
&\lesssim \sum_{k\in\mathbb Z^3}|\widehat M (t,k)\cdot \widehat{\na^{\al}{U}_0}(k)|^2\\
&\lesssim \sum_{k\in\mathbb Z^3/\{0\}}e^{-t}\Big(|\widehat{\na^{\al}{u}_0}(k)|^2+|\widehat{\na^{\al}{{\mathbb P\dive\tau}}_0}(k)|^2\Big)\\
&\lesssim e^{-t}\big(\|\na^{\al}{u}_0\|^2_{L^2(\mathbb T^3)} +\|\na^{\al+1}{\tau}_0\|^2_{L^2(\mathbb T^3)}\big),\quad \forall \al \in\N^3.
\end{split}
\eeq
Applying a similar argument as in the estimate of the derivatives of $(N \ast{U}_0, M \ast{U}_0)$, then we can obtain that\beq\label{e6-3}
\begin{split}
&\|\na^\al N(t-s,k)\ast H(s,k)\|_{L^2(\mathbb T^3)}^2\\
\lesssim &\sum_{k\in\mathbb Z^3/\{0\}}e^{-{(t-s)}}\Big(|\widehat{\na^{\al}(u\cdot \na u)}(s,k)|^2+|\widehat{\na^{\bar\al}\mathbb P\dive(u\cdot \na\tau+ (\tr\tau)\tau+Q(\tau,\na u))}(s,k)|^2\Big)\\
\lesssim &e^{-{(t-s)}}\Big(\|{\na^{\al}(u\cdot \na u)}(s)\|_{L^2}^2+\|\na^{\bar\al}\mathbb P\dive(u\cdot \na\tau+ (\tr\tau)\tau+Q(\tau,\na u))(s)\|_{L^2}^2\Big),\quad \forall \al \in\N^3,
\end{split}
\eeq
where $\bar\al \in\{\al,\al-1,\al-2\}$, and\beq\label{e6-4}
\begin{split}
&\|\na^{\al} M(t-s,k)\ast H(s,k)\|_{L^2(\mathbb T^3)}^2\\
\lesssim &\sum_{k\in\mathbb Z^3/\{0\}}e^{-{(t-s)}}\Big(|\widehat{\na^{\al}(u\cdot \na u)}(s,k)|^2+|\widehat{\na^{\al}\mathbb P\dive(u\cdot \na\tau+ (\tr\tau)\tau+Q(\tau,\na u))}(s,k)|^2\Big)\\
\lesssim &e^{-{(t-s)}}\Big(\|\na^{\al}(u\cdot \na u)(s)\|_{L^2}^2+\|\na^{\al}\mathbb P\dive(u\cdot \na\tau+ (\tr\tau)\tau+Q(\tau,\na u))(s)\|_{L^2}^2\Big),\quad \forall \al \in\N^3.
\end{split}
\eeq

Specially, together with \eqref{e6-1}, \eqref{e6-2},\eqref{e6-3} and \eqref{e6-4}, we get that\beno
\begin{split}
\|u(t)\|_{H^1}
\lesssim & e^{-\f t2} \mathcal E^\f12(0)+\sum_{\al=0}^2\int_0^t e^{-\f{t-s}2}\|\na^{\al}(u\cdot \na u)(s)\|_{L^2}ds +\int_0^t e^{-\f{t-s}2}\|\mathbb P\dive(u\cdot\na\tau)(s)\|_{L^2}ds\\
&+\int_0^t e^{-\f{t-s}2}\|\mathbb P\dive((\tr\tau)\tau)(s)\|_{L^2}ds+\int_0^t e^{-\f{t-s}2}\|\mathbb P\dive   Q(\tau,\na u)(s)\|_{L^2}ds \\
=&e^{-\f t2} \mathcal E^\f12(0)+O_1+O_2+O_3+O_4.
\end{split}
\eeno
Applying the Lemma \ref{lemm}, we verify that
\beno
\begin{split}
O_1\lesssim&\int_0^t e^{-\f{t-s}2}\|u(s)\|_{H^2}\|\na u(s)\|_{H^2}ds\\
\lesssim&\sup_{o\leq s\leq t}\|u\|_{H^2}\Big(\int_0^t(1+c_0t)^{3-\ep}\|\na^3 u(s)\|_{L^2}^2ds\Big)^\f12\Big(\int_0^t e^{-{(t-s)}}(1+c_0t)^{-3+\ep}ds\Big)^\f12\\
\lesssim&\mathcal E_1^\f12(t)\mathcal E_2^\f12(t)\Big(c_0^{-\f12} e^{-\f t4}+(1+c_0t)^{-\f32+\f\ep2}\Big).
\end{split}
\eeno
Taking advantage of the Lemma \ref{lem}, we deduce that \beno
\begin{split}
O_2\lesssim&\int_0^t e^{-\f{t-s}2}\big(\|u\|_{L^\infty}\|\na\mathbb P\dive\tau(s)\|_{L^2}+\|\na u\|_{L^6}\|\na\tau\|_{L^3}\big)ds
\lesssim \mathcal E_1^\f12(t)\mathcal E_2^\f12(t)\Big(c_0^{-\f12} e^{-\f t4}+(1+c_0t)^{-\f32+\f\ep2}\Big),
\end{split}
\eeno
 \beno
\begin{split}
O_3\lesssim&\int_0^t e^{-\f{t-s}2}\big(\|\tr\tau\|_{L^6}\|\mathbb P\dive \tau(s)\|_{L^3}+\|\na\tr\tau\|_{L^2}\|\tau(s)\|_{L^\infty}\big)ds\\
\lesssim& c_0^\f12\mathcal E_1^\f12(t)\mathcal E_4^\f12(t)\Big(c_0^{-\f12} e^{-\f t4}+(1+c_0t)^{-\f32+\f\ep2}\Big).
\end{split}
\eeno
and\beno
\begin{split}
O_4\lesssim&\int_0^t e^{-\f{t-s}2}\big(\|\na \tau\|_{L^3}\|\na u(s)\|_{L^6}+\|\tau\|_{L^\infty}\|\na^2 u(s)\|_{L^2}\big)ds
\lesssim \mathcal E_1^\f12(t)\mathcal E_2^\f12(t)\Big(c_0^{-\f12} e^{-\f t4}+(1+c_0t)^{-\f32+\f\ep2}\Big).
\end{split}
\eeno
Summing up the above estimates, we obtain that\beq\label{decay-u}
\begin{split}
\|u(t)\|_{H^2}\lesssim&
e^{-\f t2} \mathcal E^\f12(0)+\Big( \mathcal E_1^\f12(t)\mathcal E_2^\f12(t)+c_0^\f12\mathcal E_1^\f12(t)\mathcal E_4^\f12(t)\Big)\Big(c_0^{-\f12} e^{-\f t4}+(1+c_0t)^{-\f32+\f\ep2}\Big).
\end{split}
\eeq
Thus, we also get that\beno
\begin{split}
\|\mathbb P\dive\tau(t)\|_{L^2}
\lesssim & e^{-\f t2} \mathcal E^\f12(0)+\int_0^t e^{-\f{t-s}2}\|(u\cdot \na u)(s)\|_{L^2}ds +\int_0^t e^{-\f{t-s}2}\|\mathbb P\dive(u\cdot \na\tau)(s)\|_{L^2}ds\\
&+\int_0^t e^{-\f{t-s}2}\| \mathbb P\dive((\tr\tau)\tau)(s)\|_{L^2}ds+\int_0^t e^{-\f{t-s}2}\|\mathbb P\dive Q(\tau,\na u)(s)\|_{L^2}ds,
\end{split}
\eeno
which implies that\beq\label{decay-tau}
\begin{split}
\|\mathbb P\dive\tau(t)\|_{L^2}\lesssim&
e^{-\f t2} \mathcal E^\f12(0)+\Big( \mathcal E_1^\f12(t)\mathcal E_2^\f12(t)+c_0^\f12\mathcal E_1^\f12(t)\mathcal E_4^\f12(t)\Big)\Big(c_0^{-\f12} e^{-\f t4}+(1+c_0t)^{-\f32+\f\ep2}\Big).
\end{split}
\eeq

Finally, we have\beq\label{E3}
\mathcal E_3(t)\lesssim \mathcal E(0)+c_0^{-3+\ep}\mathcal E_1(t)\mathcal E_2(t)+c_0^{-2+\ep}\mathcal E_1(t)\mathcal E_4(t).
\eeq

\subsection{The estimates of $\mathcal E_4(t)$ and $\mathcal E_5(t)$}
From \eqref{tau}, we have \beno
\tr \tau(t, x)=\f{\tr \tau_0(q_t^{-1}(x))}{1+\tr\tau_0(q_t^{-1}(x))t}, \quad \forall t\in[0,T].
\eeno
which implies that\beno
\begin{split}
\na \tr \tau(t, x)=\f{\big(\na q_t^{-1}(x)\cdot\na\big)\tr \tau_0(q_t^{-1}(x))}{\big(1+\tr\tau_0(q_t^{-1}(x))t\big)^2},
\end{split}
\eeno
and\beno
\begin{split}
\na^2 \tr \tau(t, x)=&\f{\pa_k\pa_l\tr \tau_0(q_t^{-1}(x))\cdot\na (q_t^{-1}(x))^k\cdot\na(q_t^{-1}(x))^l}{\big(1+\tr\tau_0(q_t^{-1}(x))t\big)^2}+\f{\pa_k\tr \tau_0(q_t^{-1}(x))\cdot\na^2 (q_t^{-1}(x))^k}{\big(1+\tr\tau_0(q_t^{-1}(x))t\big)^2}\\
&\quad-\f{2t\cdot\pa_k\tr \tau_0(q_t^{-1}(x))\cdot\pa_l\tr \tau_0(q_t^{-1}(x))\cdot\na (q_t^{-1}(x))^k\cdot\na(q_t^{-1}(x))^l}{\big(1+\tr\tau_0(q_t^{-1}(x))t\big)^3}.
\end{split}
\eeno
Taking $L^2$-norm to the above equations and using the fact that $\dive u=0$, then we deduce that\beno
\begin{split}
\|\na \tr \tau(t, x)\|_{L^2}^2&\lesssim\|\na q_t^{-1}(x)\|_{L^\infty}^2\int_{\mathbb T^3} \f{|\na\tr \tau_0(q_t^{-1}(x))|^2}{|1+\tr\tau_0(q_t^{-1}(x))t|^4}dx\\
&\lesssim\|\na q_t^{-1}(x)\|_{L^\infty}^2\int_{\mathbb T^3} \f{|\na\tr \tau_0(x)|^2}{|1+\tr\tau_0(x)t|^4}dx\lesssim\f{\|\na\tr\tau_0\|_{L^2}^2\|\na q_t^{-1}(x)\|_{L^\infty}^2}{|1+c_0t|^4},
\end{split}
\eeno
and\beno
\begin{split}
\|\na^2 \tr \tau(t, x)\|_{L^2}^2&\lesssim\|\na q_t^{-1}(x)\|_{L^\infty}^4\int_{\mathbb T^3} \f{|\na^2\tr \tau_0(x)|^2}{|1+\tr\tau_0(x)t|^4}dx+\|\na^2 q_t^{-1}(x)\|_{L^\infty}^2\int_{\mathbb T^3} \f{|\na\tr \tau_0(x)|^2}{|1+\tr\tau_0(x)t|^4}dx\\
&\quad +4 \|\na q_t^{-1}(x)\|_{L^\infty}^4\int_{\mathbb T^3} \f{t^2|\na\tr \tau_0(x)|^4}{|1+\tr\tau_0(x)t|^6}dx\\
&\lesssim \f{\big(\|\na^2\tr\tau_0\|_{L^2}^2+c_0^{-2}\|\na\tr\tau_0\|_{L^4}^4\big)\|\na q_t^{-1}(x)\|_{L^\infty}^4}{|1+c_0t|^4}+\f{\|\na\tr\tau_0\|_{L^2}^2\|\na^2 q_t^{-1}(x)\|_{L^\infty}^2}{|1+c_0t|^4}.
\end{split}
\eeno
Applying the Lemma \ref{qt}, then we get
\beno
\begin{split}
&\quad\int_0^t (1+c_0s)^{3-\ep}\|\na \tr \tau(s)\|_{L^2}^2ds\\
&\leq \|\na\tr\tau_0\|_{L^2}^2\|\na q_t^{-1}(x)\|_{L^\infty}^2\int_0^t (1+c_0s)^{-1-\ep}ds\\
&\lesssim c_0^{-1}\|\na\tr\tau_0\|_{L^2}^2\exp\Big(\int_0^t \|\na u(s)\|_{L^\infty}ds\Big)\\
&\lesssim c_0^{-1}\|\na\tr\tau_0\|_{L^2}^2\exp\Big\{\Big(\int_0^t (1+c_0s)^{3-\ep}\|\na^3u(s)\|_{L^2}^2ds\Big)^\f12\Big(\int_0^t (1+c_0s)^{-3+\ep}ds\Big)^\f12\Big\}\\
&\lesssim c_0^{-1}\|\na\tr\tau_0\|_{L^2}^2\exp\Big(c_0^{-\f12}\mathcal E_2^{\f12}(t)\Big),
\end{split}
\eeno
which implies \beq\label{E4}
\mathcal E_4(t)\lesssim c_0^{-2}\mathcal E(0)\exp\Big(c_0^{-\f12}\mathcal E_2^{\f12}(t)\Big).
\eeq

Applying the Lemma \ref{qt} again, then we get
\beno
\begin{split}
&\quad\int_0^t (1+c_0s)^{3-\ep}\|\na^2 \tr \tau(s)\|_{L^2}^2ds\\
&\lesssim c_0^{-1}\big(\|\na^2\tr\tau_0\|_{L^2}^2+c_0^{-2}\|\na\tr\tau_0\|_{L^2}\|\na^2\tr\tau_0\|_{L^2}^3\big)\|\na q_t^{-1}(x)\|_{L^\infty}^4+c_0^{-1}\|\na\tr\tau_0\|_{L^2}^2\|\na^2 q_t^{-1}(x)\|_{L^\infty}^2\\
&\lesssim c_0^{-1}\big(\|\na^2\tr\tau_0\|_{L^2}^2+c_0^{-2}\|\na\tr\tau_0\|_{L^2}\|\na^2\tr\tau_0\|_{L^2}^3\big)\exp \Big(\int_0^t \|\na u(s)\|_{L^\infty}ds\Big)\\
&\qquad+c_0^{-1}\|\na\tr\tau_0\|^2_{L^2}\exp \Big(\int_0^t \|\na u(s)\|_{L^\infty}ds\Big)\int_0^t \|\na^2 u(s)\|_{L^\infty}\exp \Big(\int_0^s\|\na u(s')\|_{L^\infty}ds'\Big)ds\\
&\lesssim c_0^{-1}\big(\|\na^2\tr\tau_0\|_{L^2}^2+c_0^{-2}\|\na\tr\tau_0\|_{L^2}\|\na^2\tr\tau_0\|_{L^2}^3\big)\exp\Big(c_0^{-\f12}\mathcal E_2^{\f12}(t)\Big)\\
&\qquad+c_0^{-1}\|\na\tr\tau_0\|^2_{L^2}\exp\Big(c_0^{-\f12}\mathcal E_2^{\f12}(t)\Big)\int_0^t \|\na^2 u(s)\|_{L^\infty}ds.
\end{split}
\eeno
Now we estimate the $\int_0^t \|\na^2 u(s)\|_{L^\infty} ds$. For this purpose, we rewrite the first equation of the system \eqref{PTT} as following:
\beno
u_{t}-\tri u=\mathbb P\dive \tau-\mathbb P(u\cdot \na u),
\eeno
then we have\beno
u(t,x)=e^{-t\tri}u_0(x)+\int_0^t e^{-(t-s)\tri}(\mathbb P\dive \tau-\mathbb P(u\cdot \na u))(s)ds.
\eeno
By standard estimates about the Heat equation, we obtain that\beno
\begin{split}
\|\na^2 e^{-t\tri}u_0(x)\|_{L^\infty}
\lesssim& \sum_{k\in\mathbb Z^3/\{0\}}e^{-|k|^2t}|k|^2|\widehat u_0(k)|
\lesssim e^{-t}\|u_0\|_{L^2}.
\end{split}
\eeno
Using the \eqref{decay-tau}, then we have\beno
\begin{split}
&\|\int_0^t \na^2 e^{-(t-s)\tri}(\mathbb P\dive \tau)(s)ds\|_{L^\infty}\\
\lesssim&  \sum_{k\in\mathbb Z^3/\{0\}}\int_0^t e^{-|k|^2(t-s)}|k|^2|\widehat{\mathbb P\dive \tau}(s,k)| ds\\
\lesssim& \int_0^t e^{-(t-s)}\|{\mathbb P\dive \tau}(s)\|_{L^2}ds\\
\lesssim& \mathcal E^\f12(0)\int_0^t e^{-(t-s)}e^{-\f s2}ds+c_0^{-\f12}\Big(\mathcal E_1^\f12(t)\mathcal E_2^\f12(t)+c_0^{\f12}\mathcal E_1^\f12(t)\mathcal E_4^\f12(t)\Big)\int_0^te^{-(t-s)}e^\f s4 ds\\
&\quad+\Big(\mathcal E_1^\f12(t)\mathcal E_2^\f12(t)+c_0^{\f12}\mathcal E_1^\f12(t)\mathcal E_4^\f12(t)\Big)\int_{0}^te^{-(t-s)}(1+c_0s)^{-\f32+\f\ep2} ds\\
\lesssim& e^{-\f t2}\mathcal E^\f12(0)+c_0^{-1}\Big(\mathcal E_1^\f12(t)\mathcal E_2^\f12(t)+c_0^{\f12}\mathcal E_1^\f12(t)\mathcal E_4^\f12(t)\Big)e^{-\f t4}+\Big(\mathcal E_1^\f12(t)\mathcal E_2^\f12(t)+c_0^{\f12}\mathcal E_1^\f12(t)\mathcal E_4^\f12(t)\Big)(1+c_0t)^{-\f32+\f\ep2} ,
\end{split}
\eeno
and\beno
\begin{split}
\|\int_0^t \na^2 e^{-(t-s)\tri}\mathbb P(u\cdot \na u)(s)ds\|_{L^\infty}
\lesssim&  \sum_{k\in\mathbb Z^3/\{0\}}\int_0^t e^{-|k|^2(t-s)}|k|^2|\widehat{\mathbb P(u\cdot \na u)}(s,k)| ds\\
\lesssim& \int_0^t e^{-(t-s)}\|{u\cdot \na u}(s)\|_{L^2} ds\\
\lesssim& \mathcal E_1^\f12(t) \mathcal E_2^\f12(t) \Big( \int_0^te^{-2(t-s)}(1+c_0s)^{-3+\ep}ds\Big)^\f12\\
\lesssim &c_0^{-\f12}\mathcal E_1^\f12(t) \mathcal E_2^\f12(t)e^{-\f t2}+\mathcal E_1^\f12(t) \mathcal E_2^\f12(t)(1+c_0s)^{-\f32+\f\ep2}.
\end{split}
\eeno
Together with the above estimates, we have\beno
\|\na^2 u(t)\|_{L^\infty}\lesssim e^{-\f t2}\mathcal E^{\f12}(0)+c_0^{-1}\Big(\mathcal E_1^\f12(t)\mathcal E_2^\f12(t)+c_0^{\f12}\mathcal E_1^\f12(t)\mathcal E_4^\f12(t)\Big)e^{-\f t4}+\Big(\mathcal E_1^\f12(t)\mathcal E_2^\f12(t)+c_0^{\f12}\mathcal E_1^\f12(t)\mathcal E_4^\f12(t)\Big)(1+c_0t)^{-\f32+\f\ep2} ,
\eeno
integrating with time, then we get that\beq\label{na2u}
\int_0^t \|\na^2 u(s)\|_{L^\infty}ds\lesssim\mathcal E^{\f12}(0)+c_0^{-1}\mathcal E_1^\f12(t)\mathcal E_2^\f12(t)+c_0^{-\f12}\mathcal E_1^\f12(t)\mathcal E_4^\f12(t).
\eeq
Finally, we have\beno
\begin{split}
&\int_0^t (1+cs)^{3-\ep}\|\na^2 \tr \tau(s)\|_{L^2}^2ds\\
\lesssim& c_0^{-1}\big(\|\na^2\tr\tau_0\|_{L^2}^2+c_0^{-2}\|\na\tr\tau_0\|_{L^2}\|\na^2\tr\tau_0\|_{L^2}^3\big)\exp\Big(c_0^{-\f12}\mathcal E_2^{\f12}(t)\Big)\\
&\qquad+c_0^{-1}\|\na\tr\tau_0\|^2_{L^2}\exp\Big(c_0^{-\f12}\mathcal E_2^{\f12}(t)\Big)\Big(\mathcal E^{\f12}(0)+c_0^{-1}\mathcal E_1^\f12(t)\mathcal E_2^\f12(t)+c_0^{-\f12}\mathcal E_1^\f12(t)\mathcal E_4^\f12(t)\Big).
\end{split}
\eeno
Combining with the above estimates yields that\beq\label{E5}
\begin{split}
\mathcal E_{5}(t)
\lesssim& \big(\tilde{\mathcal E}(0)+c_0^{-1}\mathcal E^\f12(0)\tilde{\mathcal E}^\f32(0)\big)\exp\Big(c_0^{-\f12}\mathcal E_2^{\f12}(t)\Big)+c_0^{-2}\mathcal E(0)\exp\Big(c_0^{-\f12}\mathcal E_2^{\f12}(t)\Big)\\
&\quad\times\Big(\mathcal E^{\f12}(0)+c_0^{-1}\mathcal E_1^\f12(t)\mathcal E_2^\f12(t)+c_0^{-\f12}\mathcal E_1^\f12(t)\mathcal E_4^\f12(t)\Big).
\end{split}
\eeq

\subsection{Proof of the main theorem}
\begin{proof}
In this subsection, we will combine the above a priori estimates of $\mathcal E_1(t)$, $\mathcal E_2(t)$, $\mathcal E_3(t)$, $\mathcal E_4(t)$, $\mathcal E_5(t)$ together and give the proof of the Theorem \ref{them}. For any $t\in[0,T]$, we have
\beq\label{Etotal1}
\begin{split}
\Lambda_1(t)=&\mathcal E_1(t)+\mathcal E_2(t)\\
\leq&C^*\mathcal E(0)+Cc_0^{-\f12}\Big(\mathcal E_1^\f32(t)+\mathcal E_2^\f32(t)\Big)+Cc_0^{\f12}\Big(\mathcal E_1(t)+\mathcal E_2(t)\Big)\Big(\mathcal E_4^\f12(t)+\mathcal E_5^\f12(t)\Big)\\
&\quad+C\mathcal E_1(t)\Big(\mathcal E_4^\f14(t)\mathcal E_5^\f14(t)+\mathcal E_5^\f12(t)\Big),
\end{split}
\eeq
\beq\label{Etotal2}
\begin{split}
\Lambda_2(t)=&\mathcal E_{4}(t)\lesssim c_0^{-2}\mathcal E(0)\exp \Big(c_0^{-\f12}\mathcal E_2^{\f12}(t)\Big),
\end{split}
\eeq
\beq\label{Etotal3}
\begin{split}
\Lambda_3(t)=\mathcal E_{5}(t)
\lesssim &\big(\tilde{\mathcal E}(0)+c_0^{-1}\mathcal E^\f12(0)\tilde{\mathcal E}^\f32(0)\big)\exp\Big(c_0^{-\f12}\mathcal E_2^{\f12}(t)\Big)+c_0^{-2}\mathcal E(0)\exp\Big(c_0^{-\f12}\mathcal E_2^{\f12}(t)\Big)\\
&\quad\times\Big(\mathcal E^{\f12}(0)+c_0^{-1}\mathcal E_1^\f12(t)\mathcal E_2^\f12(t)+c_0^{-\f12}\mathcal E_1^\f12(t)\mathcal E_4^\f12(t)\Big).
\end{split}
\eeq
Due to the local existence theory, there exists a positive time $T$ such that\beq\label{lambda}
\Lambda_1(t) \leq 2C^*\d^2_0,\quad \Lambda_2(t) \leq C^*,\quad \Lambda_3(t) \leq C^*(\tilde\ep_0+\d_0^\f12),\quad \forall t\in[0,T].
\eeq
Let $T^*$ be the maximal time for what \eqref{lambda} holds.
Since $\mathcal E(0)\leq \d_0^2\leq\ep_0,$ $\tilde{\mathcal E}(0)\leq \tilde\ep_0$, $c_0=\f12\d_0$, then it follows that for any $t\in[0,T^*)$
\beno
c_0^{-\f12}\Big(\mathcal E_1^\f32(t)+\mathcal E_2^\f32(t)\Big)
\lesssim \d_0^{\f52},
\eeno
\beno
c_0^{\f12}\Big(\mathcal E_1(t)+\mathcal E_2(t)\Big)\Big(\mathcal E_4^\f12(t)+\mathcal E_5^\f12(t)\Big)
\lesssim c_0^{\f12}\Big(\mathcal E_1(t)+\mathcal E_2(t)\Big),
\eeno
\beno
\mathcal E_1(t)\Big(\mathcal E_4^\f14(t)\mathcal E_5^\f14(t)+\mathcal E_5^\f12(t)\Big)
\lesssim (\tilde\ep_0+\d_0^\f12)^\f14\mathcal E_1(t),
\eeno
\beno
\big(\tilde{\mathcal E}(0)+c_0^{-1}\mathcal E^\f12(0)\tilde{\mathcal E}^\f32(0)\big)\exp\Big(c_0^{-\f12}\mathcal E_2^{\f12}(t)\Big)
\lesssim \tilde\ep_0,
\qquad
\mathcal E^{\f12}(0)+c_0^{-1}\mathcal E_1^\f12(t)\mathcal E_2^\f12(t)+c_0^{-\f12}\mathcal E_1^\f12(t)\mathcal E_4^\f12(t) \lesssim \d_0^{\f12}.
\eeno
Under the setting of initial data, there exists  small enough numbers $\ep_0, \tilde\ep_0$ such that $\mathcal E(0)\leq \d_0^2\leq\ep_0$, $\tilde{\mathcal E}(0)\leq\tilde\ep_0$. By virtue of \eqref{Etotal1}, \eqref{Etotal2}, \eqref{Etotal3}  and the smallness assumption on $\d_0$ and $ \tilde\ep_0$, we get that\beno
\Lambda_1(t)\leq C^*\d_0^2+C\d_0^{\f52}+C\big(c_0^{\f12}+(\tilde\ep_0+\d_0^{\f12})^\f14\big)\Lambda_1(t)<2C^* \d_0^2,
\eeno
\beno
\Lambda_2(t)\leq C<C^*,
\eeno
\beno
\Lambda_3(t)\leq C\tilde\ep_0+C\d_0^{\f12}<C^*(\tilde\ep_0+\d_0^{\f12}).
\eeno
By standard continuity argument and total energy \eqref{Etotal1}, \eqref{Etotal2}, \eqref{Etotal3}, we can show that $T^*=\infty$ provided that $\d_0$ and $\tilde \ep_0$ is small enough. Moreover, we deduce the time decay estimate for the velocity $u$ and the quantity $P\dive\tau$ from \eqref{E3} that \beno
\begin{split}
\|u\|_{H^2}+\|\mathbb P\dive\tau\|_{L^2}\lesssim(1+t)^{-\f32+\f\ep2}.
\end{split}
\eeno
Hence, we finish the proof of the Theorem \ref{them}.

\end{proof}

{\bf Acknowledgements}.
 This work was partially supported by NSF of China (Grants No.11701586, No.11671407, No.11431015 and No.11801586) and the Central Universities (Grants No. 18lgpy66).

\bibliographystyle{abbrv} 
\bibliography{PTTref}

\end{document}